\documentclass[10pt]{amsart}

\usepackage[all, cmtip]{xy}
\usepackage{amsmath,amsthm,amssymb,amsfonts,amscd}
\usepackage{tikz-cd}
\usepackage{float} 

\newtheorem{theorem}{Theorem}[section]
\newtheorem{lemma}[theorem]{Lemma}
\newtheorem{proposition}[theorem]{Proposition}
\newtheorem{corollary}[theorem]{Corollary}

\theoremstyle{definition}
\newtheorem{definition}[theorem]{Definition}
\newtheorem{notation}[theorem]{Notation}
\newtheorem{example}[theorem]{Example}

\newtheorem{question}[theorem]{Question}

\theoremstyle{remark}
\newtheorem{remark}[theorem]{Remark}

\numberwithin{equation}{section}

\begin{document}

\title{Symmetric and K\"ahler--Einstein Fano polygons}
 
\author{DongSeon Hwang}
\address{Center for Complex Geometry, Institute for Basic Science (IBS), Daejeon 34126, Republic of Korea}
\email{dshwang@ibs.re.kr}

\author{Yeonsu Kim}
\address{Department of Mathematics, Ajou University, Suwon 16499, Republic Of Korea} 
\email{kys3129@ajou.ac.kr} 
 
\subjclass[2020]{Primary 52B20; Secondary 52B15,  32Q20, 14M25,  14J45.}
 
\date{April 2, 2024 and, in revised form, Sep 4, 2024}

\keywords{Fano polytope, symmetric polytope, automorphism, toric del Pezzo surface, K\"ahler--Einstein metric} 

\begin{abstract}
We investigate \emph{singular} symmetric and K\"ahler--Einstein Fano polytopes. More precisely, we show that every symmetric Fano polytope is K\"ahler--Einstein generalizing the work by Batyrev and Selivanova, and study the automorphism groups of symmetric and K\"ahler--Einstein Fano polygons in detail. In particular, every finte subgroup of $GL_2(\mathbb{Z})$ is an automorphism group of a K\"ahler--Einstein Fano polygon.
\end{abstract}

\maketitle
 
\section {Introduction}
It has been an important problem to determine which Fano variety admits a K\"ahler-Einstein metric. This differential geometric problem  can  now be understood purely in terms of algebraic geometry. More precisely, the existence of the K\"ahler-Einstein metric for a Fano variety $X$ turned out to be equivalent to the K-polystability of $X$ (\cite{CDS15a}, \cite{CDS15b}, \cite{CDS15c}, \cite{Tian15}, \cite{B} and \cite{LXZ}).  

In this note we consider toric Fano varieties. Recall that the class of toric Fano varieties  up to isomorphism has a one-to-one correspondence with the class of Fano polytopes  up to unimodular transformation. Here, a {\em Fano polytope} of dimension $n$ is a full dimensional convex lattice polytope in $\mathbb{R}^n$  such that the vertices are primitive lattice points and the origin is an interior point. See \cite{KN}  for more details. Under this one-to-one correspondence, many of the algebro-geometric properties can be understood by means of convex geometry. Now it is natural to ask the following question.
 
\begin{question}\label{CG=>KE}
    Which convex geometric properties of a Fano polytope determine the existence of a K\"ahler--Einstein metric on the corresponding toric Fano variety?
\end{question}

One complete  answer to Question \ref{CG=>KE} is given by the   following purely convex geometric characterization for the existence of the K\"ahler-Einstein metric, which is further generalized for toric log Fano pairs in \cite[Theorem 1.2]{BB}.

\begin{theorem}$($\cite[Theorem 1.4 and Proposition 3.2]{SZ}, \cite[Theorem 1.2]{BB}\label{SZ-theorem}$)$ 
Let $X$ be a toric Fano variety. Then $X$ admits a K\"ahler--Einstein metric if and only if the barycenter of its moment polytope   is the origin.    
\end{theorem} 

On the other hand, the existence of  a K\"ahler--Einstein metric is governed by certain properties of the automorphism group. For example,  the automorphism group of Fano variety admitting a K\"ahler--Einstein metric is reductive (\cite[Th\'eor\`eme 1]{Matsushima} and \cite[Theorem 1.3]{ABHLX}). So one might consider a variation of  Question \ref{CG=>KE} in terms of the automorphism group.

\begin{question}\label{Aut=>KE}
    Which properties of the automorphism group of a Fano polytope determine the existence of a K\"ahler--Einstein metric on the corresponding toric Fano variety?
\end{question}

A theorem of Batyrev and Selivanova can be regarded as an answer for this question. We first recall some terminologies. A Fano polytope $P$, or the corresponding toric Fano variety,  is said to be {\em symmetric} if the origin is the only  fixed point under the automorphism group $Aut(P)$  of $P$. A Fano polytope $P \subset N_\mathbb{R}$ is said to be {\em smooth} if the vertices of any facet of $P$ form a $\mathbb{Z}$-basis for $N$, and {\em singular} if otherwise. Motivated by Theorem \ref{SZ-theorem}, a Fano polytope is said to be {\em K\"ahler--Einstein} if the barycenter of its moment polytope is the origin. From now on, we will use terminologies in convex geometry.
 
\begin{theorem}\cite[Theorem 1.1]{BS}\label{BS-Theorem}
Let $X$ be a smooth Fano polytope. If $X$ is symmetric, then $X$ is 
  K\"ahler-Einstein. 
\end{theorem}

There had been some interests in the converse direction of Theorem \ref{BS-Theorem}. It turned out that, for smooth Fano polytopes, the class of K\"ahler-Einstein Fano polytopes coincides with that of symmetric Fano polytopes  up to dimension six by \cite[Proposition 2.1]{NP}, and, for each $n \geq 7$,  there exists an $n$-dimensional non-symmetric K\"ahler-Einstein  toric Fano manifold, which is not isomorphic to a product of lower-dimensional toric manifolds by \cite[Proposition 2.1]{NP} and \cite[Corollary 5.3]{N}.

In this note, we revisit the above problems for possibly \emph{singular} Fano polytopes. Thanks to Theorem \ref{SZ-theorem} it is now easy to prove Theorem \ref{BS-Theorem} without the smoothness assumption.

\begin{theorem}\label{symm=>KE}
    Every symmetric Fano polytope  is K\"ahler--Einstein.
\end{theorem}

Even though the converse of Theorem \ref{symm=>KE} does not hold, one may still expect that the  K\"ahler--Einstein Fano polytope should have many "symmetries". For this purpose,  we systematically study the automorphism groups of symmetric  and K\"ahler--Einstein Fano polygons.

One must be careful when considering singular symmetric Fano polytopes since their automorphism group  can have a fixed point that is not a lattice point. To clarify   we introduce the following notion.   A Fano polytope $P$ is said to be {\em lattice symmetric} if the origin is the only {\em lattice} point fixed by the automorphism group of $P$.
Note that if a Fano polytope $P$ is smooth (or reflexive in general), then $P$ is symmetric if and only if it is lattice symmetric since the only interior lattice point of $P$ is the origin by \cite[Proposition 1.12]{Nill2005}.  We start by characterizing lattice symmetric but not symmetric Fano polygons. 
 
\begin{theorem}\label{Aut(L-symm)}
Let $P$ be a lattice symmetric Fano polygon. 
Then $P$ is not symmetric if and only if $Aut(P) \cong \mathbb{Z}/2$ is generated by a reflection.  In this case  $P$ is unimodularly equivalent to the Fano polygon   $$P_{m,n}=conv\{(m+1,-m),(-m,m+1),(-n-1,n),(n,-n-1)\}$$ 
    for some integers $m$ and $n$ satisfying $m\neq n$, and it is not K\"ahler--Einstein.
\end{theorem}
 
It turned out that the automorphism group of a K\"ahler--Einstein polygon can be arbitrary. 

\begin{theorem}\label{KE_aut}
    Every finite subgroup of $GL_2(\mathbb{Z})$ is an automorphism group of a K\"ahler--Einstein Fano polygon.
\end{theorem}

In particular, one can find examples of K\"ahler--Einstein Fano polygons with trivial automorphism group (Example \ref{KEnotsymm} and \ref{Ex4}). Hence the converse of Theorem \ref{symm=>KE} does not hold even in  dimension $2$ when $P$ is singular. 

\begin{question}
    Can  Theorem \ref{KE_aut} be generalized in higher dimension?
\end{question}
   
In a previous version of this note, due to lack of examples, we raised a question asking whether every K\"ahler--Einstein Fano polygon that is not symmetric is a triangle. The authors  are informed by Thomas Hall that he has found an example of a non-symmetric  K\"ahler--Einstein Fano quadrilateral, hence the answer is negative (\cite[Proposition 1.5]{Hall}). 
  
\section{Preliminaries}
In this section we recall some basic notions of Fano  polygons and review some basic properties of   K\"ahler-Einstein  Fano polygons.
\subsection{Fano polygons}  
Let $P$ be a convex lattice polygon.  Then the dual polygon $P^\ast$ of $P$ is defined by $$P^\ast=\{w \in \mathbb{R}^2  \ | \  \langle w,v_i \rangle \geq -1 \  \textrm{for every } v_i \in P\}.$$
 
We denote by $vert(P)$ the set of all vertices of $P$. The {\em order} of the two dimensional cone $C$ spanned by two lattice points $v_i = (x_i, y_i)$ and $v_{i+1} = (x_{i+1}, y_{i+1})$, denoted by $ord(v_i,v_{i+1})$, is defined by 
$$ord(v_i,v_{i+1}):=det\left( \begin{array}{ccr} x_i & x_{i+1} \\
y_i & y_{i+1} \\ \end{array} \right)=x_iy_{i+1}-y_ix_{i+1}.$$
 
\begin{notation}
For a lattice point $v = (x, y)$, we define the {\em primitive index} $I(v)$ of $v$ by $I(v)=gcd(x,y)$.
\end{notation}
A lattice point  is said to be {\em primitive} if its primitive index  is one. A convex lattice polytope $P$ is said to be  {\em Fano} if the vertices of $P$ are primitive and the origin is an interior point of $P$. A $\mathbb{Z}$-linear tranformation is said to be {\em unimodular} if the corresponding matrix has determinant $\pm 1$.  

\begin{lemma}\label{primitivity}
The primitive index of a lattice point is invariant under a unimodular transformation. 
\end{lemma}

\begin{proof}  
Let $(x,y) $ be a lattice point and $(x', y')$ be its image under a unimodular transformation. Then we can write 
$x' = ax+by$ and $y'=cx+dy$ for some integers $a,b,c,d$ satisfying $|ad-bc| = 1$. Then $gcd(x,y)$ divides $gcd(x',y')$. Since any unimodular transformation is invertible,   $gcd(x',y')$ divides $gcd(x,y)$. 
\end{proof}

Every affine toric surface is either an affine plane $\mathbb{A}^2$ or a cyclic quotient singularity. This can be easily derived   from the following   lemma below, which will be used later.

\begin{lemma}\label{fix_pt}\cite[Proposition 10.1.1]{CLS}
For any primitive points $v_1$ and $v_2$, we may take $v_1=(0,1)$ and $v_2=(n,-k)$ for some integers $n$ and $k$ satisfying $0 < k \leq n$ and $gcd(n,k)=1$ up to  unimodular transformation.
\end{lemma}

\begin{notation}
    \begin{enumerate}
        \item The affine cone generated by two vectors $v$ and $w$ in $\mathbb{R}^2$ is denoted by $cone(v,w)$.
        \item For relatively prime integers $k$ and $n$ with  $0 < k \leq n$, we denote by $\frac{1}{n}(1,k)$ the affine cone $cone((0,1),(n,-k))$. In particular,  $\frac{1}{n}(1,n-1)$ is also denoted by $A_n$.
    \end{enumerate}
\end{notation}

 \begin{remark}
     The affine toric surface corresponding to the affine cone $\frac{1}{n}(1,k)$ is the quotient of $\mathbb{A}^2$ by the action $(x,y) \mapsto (\xi x, \xi^k y)$ where $\xi$ is a primitive $n$-th root of unity. It has a unique quotient singular point at the origin of type $\frac{1}{n}(1,k)$.
 \end{remark}

Unless otherwise stated, we always list vertices of a lattice polygon in a counterclockwise order.

\subsection{K\"ahler-Einstein Fano polygons}
One can compute the barycenter of a convex lattice polygon using the following lemma.
\begin{lemma}\label{Bell}
Let $P$ be a convex polygon with $n$ vertices $v_1$, \ldots, $v_n$ written in counterclockwise order containing the origin as an interior point. Then we have the following.
\begin{enumerate}
    \item The barycenter of $P$ is given by $$\frac{\sum\limits_{i=1}^{n}(v_i+v_{i+1})ord(v_i,v_{i+1})}{3\sum\limits_{i=1}^{n}ord(v_i,v_{i+1})}$$
    where $v_i$ is adjacent to $v_{i+1}$ for each $i$.
\item If $P$ is a triangle with vertices  $(x_1, y_1), (x_2, y_2), (x_3, y_3)$, then the barycenter of $P$  is given by $(\frac{x_1+x_2+x_3}{3}, \frac{y_1+y_2+y_3}{3})$.  
\end{enumerate}
\end{lemma}

\begin{proof}
    This is an easy calculus exercise. (1) Let $T_{i}$ be the triangle  with vertices $v_i$, $v_j$, and the origin. Write $v_i =(x_i, y_i)$ for each $i$. Then it is easy to compute that the barycenter of $T_{i}$ is given by $(\frac{x_i+x_{i+1}}{3}, \frac{y_i+y_{i+1}}{3}).$ Denote by $A(T_i)$ the area of $T_i$ for each $i$ and $A(P)$ the area of $P$. Now let $(b_x, b_y)$ be the barycenter of $P$. Then $b_x$ is given by 
    $$b_x = \frac{1}{A(P)}\underset{i=1}{\overset{n}{\sum}} \frac{(x_i+x_{i+1})}{3}A(T_i) = \frac{\underset{i=1}{\overset{n}{\sum}} (x_1+x_{i+1})ord(v_i,v_{i+1})}{3\underset{i=1}{\overset{n}{\sum}} ord(v_i, v_{i+1})}$$ 
    where $x_{n+1} = x_1$ and $y_{n+1} = y_1$. Similarly we get $b_y = \frac{\underset{i=1}{\overset{n}{\sum}} (y_1+y_{i+1})ord(v_i,v_{i+1})}{3\underset{i=1}{\overset{n}{\sum}} ord(v_i, v_{i+1})}$.\\
    (2) directly follows from (1).
\end{proof}

Since $ord(v_i,v_{i+1})$ is invariant under a unimodular transformation for every $i$, we immediately get the following.

\begin{corollary}\label{multi}
The vanishing of the barycenter of a lattice polygon is invariant under a unimodular transformation.  
\end{corollary}
  
Now one can check the property of being K\"ahler-Einstein only by looking at the vertices of $P$.

\begin{lemma}\label{dual_vert} Let $P$ be a Fano polygon with vertices $v_1, \ldots, v_n$ written in counterclockwise order where $v_i=(x_i,y_i)$ for $i=1, \dots ,n$. Let   $a_{i,i+1}:=ord(v_i,v_{i+1})$. Then we have the following.
\begin{enumerate}
    \item Each vertex $v^\ast_i$ of the dual polygon $P^\ast$ of $P$ is given by $$v^\ast_i=\frac{1}{a_{i,i+1}}(y_{i}-y_{i+1},x_{i+1}-x_{i}).$$
    \item $P$ is  K\"ahler-Einstein if and only if $$\sum\limits_{i=1}^{n}ord(v^\ast_i,v^\ast_{i+1})(v^\ast_i+v^\ast_{i+1})=(0,0).$$ 
  if and only if
$$ \sum\limits_{i=1}^{n} \frac{a_{i,i+1}+a_{i+1,i+2}+a_{i+2,i}}{a_{i,i+1}a_{i+1,i+2}}\Big(\frac{1}{a_{i,i+1}}(v_i-v_{i+1})+\frac{1}{a_{i+1,i+2}}(v_{i+1}-v_{i+2})\Big)=(0,0).$$
\end{enumerate}
\end{lemma}

\begin{proof} 
Direct calculation shows (1).  The first equivalence of (2) follows from Theorem \ref{SZ-theorem} and Corollary \ref{multi}. By using (1),  it is not hard to see the last equivalence.  
\end{proof}

\section{Symmetric and K\"ahler-Einstein Fano polygons}

Let $N$ be a lattice and $P$ be a Fano polytope in $N_\mathbb{R}:=N \otimes_\mathbb{Z} \mathbb{R}$. An \emph{automorphism} of   $P$   is a $GL_{\mathbb{Z}}(N)$-symmetry preserving $P$. Let $M$ be the dual lattice of $N$ and $P^\ast$ be the dual polytope of $P$ in $M_\mathbb{R} := M \otimes_\mathbb{Z} \mathbb{R}$.

\begin{definition}
The \emph{index} of a Fano polytope $P$ is the smallest positive integer $l$ such that $lP^*$ is a lattice polytope. In this case the corresponding toric Fano variety has \emph{Cartier index} $l$.  
\end{definition}

We first observe the invariance of the symmetricity  under  dualizing for every Fano polytope as a generalization of the results of Batyrev and Nill. In fact, it is mostly straightforward to follow the argument in \cite[Proposition 5.4.2]{Nill2006} or \cite[Lemma 3.5]{NP}.

\begin{proposition}\label{symmetric-dual}
    Let $P$ be a Fano polytope of index $l$ and $P^\ast$ be its dual polytope. Then $P$ is symmetric if and only if the Fano polytope $lP^\ast$ is symmetric.
\end{proposition}

\begin{proof}
Since a lattice automorphism induces a real vectorspace automorphism which behaves well with respect to the dualizing operation from $P$ to $P^\ast$, we see that  
$$Aut_M(l P^\ast) \cong Aut_{\frac{1}{l}M}(P^\ast)  \cong  Aut_N(P)^*.$$  

Now we basically follow the proof in \cite[Proposition 5.4.2]{Nill2006} or in \cite[Lemma 3.5]{NP}. By Maschke's theorem, $N_\mathbb{R} \cong \text{Fix}_{\text{Aut}_N(P)} \oplus U$ for some $\text{Aut}_N(P)$-invariant subspace $U \subset N_\mathbb{R}$ where $\text{Fix}_{\text{Aut}_N(P)}$ denotes the set of fixed points of $N_\mathbb{R}$ under the action of $\text{Aut}_N(P)$. By dualizing we get $M_\mathbb{R} \cong \text{Fix}_{\text{Aut}_M(lP^\ast)} \oplus U^\ast$ 
where $\text{Fix}_{\text{Aut}_M(lP^\ast)}$ denotes the set of fixed points of $M_\mathbb{R}$ under the action of $\text{Aut}_M(lP^\ast)$. Thus 
 $\text{dim}_\mathbb{R} \text{Fix}_{\text{Aut}_N(P)} \leq \text{dim}_\mathbb{R} \text{Fix}_{\text{Aut}_M(lP^\ast)}.$ By symmetry, we see that $$\text{dim}_\mathbb{R} \text{Fix}_{\text{Aut}_N(P)} = \text{dim}_\mathbb{R} \text{Fix}_{\text{Aut}_M(lP^\ast)},$$
from which the result follows.
\end{proof}
 
\begin{remark} 
Proposition \ref{symmetric-dual} does not hold for lattice symmetric Fano polytopes. Take a lattice symmetric Fano polygon $P=conv\{(2,-1),(-1,2),(-1,0),(0,-1)\}.$  Then $P^\ast \cong conv\{(-1,-1),(0,1),(1,0),(1,1)\}$ is  not lattice symmetric since  $(1,1)$ and $(-1,-1)$ are fixed points of a reflection $\left( \begin{array}{ccr} 0 & 1 \\ 1 & 0 \\ \end{array} \right)$.
\end{remark}
 
Now it is easy to generalize Theorem \ref{BS-Theorem} in the singular case thanks to Theorem \ref{SZ-theorem}. 

\begin{theorem}\label{symm=>KE2}
Every symmetric Fano polytope is K\"ahler--Einstein.  
\end{theorem}

\begin{proof} 
Let $P^\ast$ be the dual polytope of the given Fano polytope of index $l$. By Proposition \ref{symmetric-dual},  $lP^\ast$ is symmetric. Thus the origin is the unique fixed point of $lP^\ast$ under the action of $Aut_M(lP^\ast)$.   Since every automorphism preserves the barycenter of a convex lattice polytope, it follows that the origin is the barycenter of $lP^\ast$. Thus the barycenter of $P^\ast$ is the origin. Now Theorem \ref{SZ-theorem} concludes the proof. 
\end{proof}
 
\section{Automorphisms}

\subsection{Automorphisms of Lattice symmetric Fano polygons} In this subsection we prove Theorem \ref{Aut(L-symm)}. We start with a  special class of  Fano polygons. 
\begin{example}[Fano polygons $P_{m,n}$]\label{SmnEX}
For non-negative integers $m$ and $n$, let $P_{m,n}:=conv\{v_1,v_2, v_3, v_4\}$ where $v_1= (m+1,-m)$, $v_2=(-m,m+1),$ $v_3=(-n-1,n)$ and $v_4=(n,-n-1).$ See Figure \ref{Smn}. For convenience we assume that $m \geq n$. It is easy to see that 
$C_1:=cone(v_1,v_2) \cong A_{2m}$, $C_3:=cone(v_3,v_4) \cong A_{2n}$,  $C_2:=cone(v_2,v_3) \cong C_4:=cone(v_4,v_1) \cong \frac{1}{m+n+1}(1,1)$ where $cone(v,w)$ denotes the cone generated by two vectors $v$ and $w$.
Then  $P_{m,n}$  is a Fano polygon corresponding to a toric del Pezzo surface of Picard number two   with at most $4$ singular points of type $$A_{2m},\frac{1}{m+n+1}(1,1),A_{2n},\frac{1}{m+n+1}(1,1).$$ 
In particular, it is reflexive, i.e., of index $1$, if $m + n \leq 1$. In other words, only $P_{0,0}$ and $P_{1,0}$ are reflexive. 
Here, $P_{0,0}$ corresponds to $\mathbb{P}^1 \times \mathbb{P}^1$  and $P_{1,0}$ corresponds to the toric del Pezzo surface of Picard number two with $3$ singular points of type $2A_1+A_2$.
\end{example}

\begin{figure}[h]
    \centering
\tikzset{every picture/.style={line width=0.75pt}} 

\begin{tikzpicture}[x=0.75pt,y=0.75pt,yscale=-1,xscale=1]

\draw  (67.5,130.86) -- (289.5,130.86)(179.61,30) -- (179.61,229.9) (282.5,125.86) -- (289.5,130.86) -- (282.5,135.86) (174.61,37) -- (179.61,30) -- (184.61,37)  ;
\draw    (79.5,70.9) -- (238.5,229.9) ;
\draw    (148.5,60.9) -- (248.5,159.9) ;
\draw    (79.5,70.9) -- (148.5,60.9) ;
\draw    (238.5,229.9) -- (248.5,159.9) ;

\draw (106,45.9) node [anchor=north west][inner sep=0.75pt]  [font=\footnotesize] [align=left] {$\displaystyle ( -m,m+1)$};
\draw (253,153.9) node [anchor=north west][inner sep=0.75pt]  [font=\footnotesize] [align=left] {$\displaystyle ( m+1,-m)$};
\draw (239,224.9) node [anchor=north west][inner sep=0.75pt]  [font=\footnotesize] [align=left] {$\displaystyle ( n,-n-1)$};
\draw (24,53.9) node [anchor=north west][inner sep=0.75pt]  [font=\footnotesize] [align=left] {$\displaystyle ( -n-1,n)$};
\end{tikzpicture}
    \caption{$P_{m,n}$}
    \label{Smn}
\end{figure}

We first compute the automorphism groups of $P_{m,n}$.

\begin{proposition}\label{Aut-P_m.n}
    The automorphism group $Aut(P_{m,n})$ is as follows.
   \begin{equation*}
Aut(P_{m,n}) \cong
\begin{cases}
\mathbb{Z}/2\mathbb{Z}, \text{ generated by a reflection,} & \text{if } m \neq n,\\
\mathbb{Z}/2\mathbb{Z} \times \mathbb{Z}/2\mathbb{Z} & \text{if } m=n \geq 1,\\
D_4 & \text{if } m=n=0.
\end{cases}
\end{equation*}
\end{proposition}

\begin{proof}
Consider the case $m \neq n$. By looking at the four cones generated by the four vertices of $P_{m,n}$, we see that the only possible nontrivial automorphism is the reflection interchanging the cones $C_2$ and $C_4$. The matrix representation of the reflection is $\left( \begin{array}{ccr} 0 & 1 \\ 1 & 0 \\ \end{array} \right)$.

Consider the case $m=n$. If $m=0$, then it is easy to see that $Aut(P_{0,0}) \cong D_4$. Assume that $m \geq 1$. Again, by looking at the four cones generated by the four vertices of $P_{m,n}$, we see that, in addition to the reflection considered above, the only possible automorphism is the rotation interchanging $C_1$ and $C_3$, and hence $C_2$ and $C_4$. The matrix representation of the rotation is $\left( \begin{array}{ccr} -1 & 0 \\ 0 & -1 \\ \end{array} \right)$.  
\end{proof} 

Now we see that $P_{m,n}$ is lattice symmetric but not K\"ahler--Einstein if $m \neq n$.

\begin{corollary}\label{m=n}
    The Fano polygon  $P_{m,n}$ is lattice symmetric. Moreover,  it is symmetric if and only if it is K\"ahler--Einstein if and only if $m = n$. 
\end{corollary}

\begin{proof}
Since a Fano polygon is symmetric if and only if it admits a non-trivial rotation,  
it is sufficient to prove that $P_{m,n}$ is K\"ahler--Einstein if and only if $m=n.$ By Lemma \ref{dual_vert},   
the dual polygon $P^\ast_{m,n}$ of $P_{m,n}$ is generated by the vertices $$\Big\{(-1,-1),\Big(\frac{m-n+1}{m+n+1},\frac{m-n-1}{m+n+1}\Big),(1,1),\Big(\frac{m-n-1}{m+n+1},\frac{m-n+1}{m+n+1}\Big)\Big\}.$$
and  the barycenter of   $P^\ast_{m,n}$ is $$\Big(\frac{m-n}{3(m+n+1)},\frac{m-n}{3(m+n+1)}\Big),$$ 
which is  the origin if and only if $m = n$. 
\end{proof}

Now we show that in fact $P_{m,n}$'s are the only lattice symmetric Fano polygons that are not symmetric.

\begin{theorem}\label{no-K.E.}
Let $P$ be a lattice symmetric Fano polygon. If $P$ is not symmetric,  $P \cong P_{m,n}$ for some non-negative integers $m$ and $n$ with $m \neq n$ up to  unimodular transformation.
\end{theorem}

\begin{proof}
Since $P$ is lattice symmetric but not symmetric, $P$ admits a non-trivial reflection. Write  $\sigma=\left( \begin{array}{ccr} 0 & 1 \\ 1 & 0 \\ \end{array} \right)$ for the matrix representation of the reflection. 
We claim that, if $(a,b)$ is a vertex of $P$, then $|a+b| = 1$.
Assume first that  $|a+b|>1$. Since $P$ admits a reflection $\sigma$, $(b,a)$ is also a vertex of $P$. Then the two dimensional cone generated by $(a,b)$ and $(b,a)$ contains a lattice point  $(1,1)$ or $(-1,-1)$, each of which is a fixed point of $\sigma$, a contradiction since $P$ is lattice symmetric. If $|a+b| = 0$, we may assume that $(a,b) = (-1, 1)$ by the  primitivity of $(a,b)$ up to the reflection $\sigma$. Then it is easy to see that    $$vert(P)=\{(0,1),(-1,1),(-1,0),(0,-1),(1,-1),(1,0)\}.$$
In this case, $P$ admits another  reflection $\left( \begin{array}{ccr} 0 & -1 \\ -1 & 0 \\ \end{array} \right)$, a contradiction since $P$ is not symmetric. This proves the claim.

Assuming the claim, by the primitivity of the vertices, we may write   $$vert(P)=\{(m+1,-m),(-m,m+1),(-n-1,n),(n,-n-1)\}$$
for some non-negative integers $m$ and $n$.
Thus $P \cong P_{m,n}$. By Corollary \ref{m=n}, we have $m \neq n$.  
\end{proof}

\subsection{Automorphisms of K\"ahler-Einstein Fano polygons}\label{KEtriangle}
We first study the automorphism groups of K\"ahler--Einstein Fano triangles. 

\begin{proposition}\label{KE-triangle}
Let $P$ be a K\"ahler-Einstein Fano triangle.
Then up to  unimodular transformation,  $$P=conv\{(a,-b),(0,1),(-a,b-1)\}$$ where $a$ and $b$ are positive integers satisfying $gcd(a,b)=gcd(a,b-1)=1$.
\end{proposition}

\begin{proof}
By Lemma \ref{fix_pt}, we may assume that $v_1=(a,-b)$ and $v_2=(0,1)$ for some positive integers $a$ and $b$ satisfying $gcd(a,b)=1$ and $a \geq b$.
Let $P=conv\{(a,-b),(0,1),(-c,d)\}$ be a K\"ahler-Einstein Fano triangle for some   integers $c$ and $d$ with $gcd(c,d)=1$. We may assume that $((a,-b),(0,1),(-c,d))$ is in a counterclockwise order. Then $ 0 < ord(v_3, v_1) = bc-ad$.
By Lemma \ref{dual_vert},  the dual polygon $P^\ast$ of $P$ can be written as  $$P^\ast = conv\Big\{\Big(-\frac{b+1}{a},-1\Big),\Big(\frac{1-d}{c},-1\Big),\Big(\frac{b+d}{bc-ad},\frac{a+c}{bc-ad}\Big)\Big\}.$$
Since $P$ is K\"ahler-Einstein, by Lemma \ref{Bell}, we see that either $c=a$ and $d = b-1$ or $c=-a$ and $d=-b$. The latter case cannot happen  since  $P$ contains the origin as an interior point.
\end{proof}

\begin{corollary}\label{Fanotriangleindex} A K\"ahler--Einstein toric del Pezzo surface $S$ of Picard number one has $3$ singular points of type 
 $\frac{1}{a}(1,b)$, $\frac{1}{a}(1,b-1)$  and $\frac{1}{a}(1,x+1)$ where $x$ is a positive integer satisfying $ax-by=1$ for some  integer $y$. 
In particular, $S$ has Cartier index $a$. 
\end{corollary}

\begin{proof}
Let $P$ be the corresponding  K\"ahler-Einstein Fano triangle. Then by Proposition \ref{KEtriangle}, 
 $P$ is generated by primitive lattice points $v_1=(a,-b)$, $v_2=(0,1)$ and $v_3=(-a,b-1)$. Since $gcd(a,b)=1$, there exist integers  $x$ and $y$ such that $-ay+bx=1$.

It is easy to see that the cone  generated by $v_1$ and $v_2$ corresponds to the singular point of type $\frac{1}{a}(1,b)$.   Next, by considering the unimodular transformation $\left( \begin{array}{ccr} -1 & 0 \\ 1 & 1 \\ \end{array} \right)$, we can easily see that the cone  generated by $v_1$ and $v_2$ corresponds to the singular point of type $\frac{1}{a}(1,b-1)$. Finally, by considering the  unimodular transformation $\left( \begin{array}{ccr} -b & -a \\ b-y_1 & a-x_1 \\ \end{array} \right)$, we see that the cone generated by $v_3$ and $v_1$ corresponds to the singular point of type $\frac{1}{a}(1,-a+x+1)$.  
\end{proof}
 
\begin{corollary}\label{KE_tri_index}
Let $P$ be a K\"ahler-Einstein Fano triangle.  Then the index of $P$ is an odd integer.
\end{corollary}

\begin{proof}
Let $P=conv\{(0,1),(a,-b),(-a,b-1)\}$ where  $a$ and $b$ are positive integers satisfying $gcd(a,b)=gcd(a,b-1)=1$ by Proposition \ref{KE-triangle}. 
By Corollary \ref{Fanotriangleindex}, the index of $P$ is $a$. If the index $a$ is an even integer, then neither $gcd(a,b)$ nor $gcd(a,b-1)$ is $1$, a contradiction. 
\end{proof}

\begin{theorem}\label{tri-auto}
Let $P$ be a K\"ahler-Einstein Fano triangle.
Then the automorphism group of $P$ is either trivial, $\mathbb{Z}_2$, $\mathbb{Z}_3$, or  $D_3$.
Moreover, the last two cases happen if and only if $P$ is   symmetric. 
\end{theorem}

\begin{proof}
The first assertion is trivial. Note that $P$ is  symmetric   if and only if $P$ admits a rotation if and only if  the automorphism group of $P$ is $\mathbb{Z}_3$ or $D_3$.  
\end{proof}

\begin{remark}\label{stdmatrix}
A K\"ahler-Einstein Fano triangle $P$ can be written as 
$$P=conv\{(a,-b),(0,1),(-a,b-1)\}$$
where $gcd(a,b)=gcd(a,b-1)=1$  by Proposition \ref{KE-triangle}.   
Then  
\begin{enumerate}
    \item The reflection fixing $v_1$ corresponds to the matrix  $\left( \begin{array}{ccr} -b+1 & -a \\ \frac{b^2-2b}{a} & b-1 \\ \end{array} \right)$. 
    \item The reflection fixing $v_2$ corresponds to the matrix $\left( \begin{array}{ccr} -1 & 0 \\ \frac{2b-1}{a} & 1 \\ \end{array} \right)$.
    \item The matrix reflection fixing $v_3$ corresponds to the matrix  $\left( \begin{array}{ccr} b & a \\ \frac{1-b^2}{a} & -b \\ \end{array} \right)$.
    \item A rotation corresponds to the matrix $\left( \begin{array}{ccr} -b & -a \\ \frac{b^2-b+1}{a} & b-1 \\ \end{array} \right)$. 
\end{enumerate} 
\end{remark}

\begin{corollary}
Let $P$ be a Fano triangle. Then  $Aut(P)\cong D_3$ if and only if 
$P=conv\{(1,0),(0,1),(-1,-1)\}$ or $P=conv\{(3,-2),(0,1),(-3,1)\}.$ The corresponding toric del Pezzo surface is $\mathbb{P}^2$ in the first case  and a cubic surface   with three $A_2$ singularities in the last case.
\end{corollary}
 
\begin{proof} It is easy to see the 'if' direction. 
Suppose that $Aut(P) \cong D_3$. We fix the vertices of $P$ as in Theorem \ref{KE-triangle}. We may assume that $a \geq b \geq 1$. Since $P$ admits a rotation, by Theorem \ref{tri-auto} (4), we see that   $a$ divides $b^2-b+1$. If $b=1$, then $a=1$ since $(a,b-1)$ is primitive. Thus, $P=conv\{(1,0),(0,1),(-1,-1)\}$. If $b=2$, then 
 $b^2-b+1=3$, so $a=3$. This leads to $P=conv\{(3,-2),(0,1),(-3,1)\}$.
Now suppose that $b\geq 3$. Since $P$ admits a reflection fixing $v_1$, by Theorem \ref{tri-auto} (1), we see that   $a$ divides $b^2-2b$, which contradicts to  $gcd(a,b)=1$.  
\end{proof}

\begin{example}\label{KEnotsymm}
    \begin{enumerate}
        \item Take $P = conv\{ (7,-3), (0,1),  (-7,2) \}$. Then $P$ is K\"ahler-Einstein and $Aut(P) \cong \mathbb{Z}/3\mathbb{Z}.$ It has $3$ singular points all of type $\frac{1}{7}(1,5).$
        \item Take $P_1 = conv\{ (5,-2), (0,1),  (-5,1) \}$. Then $P_1$ is K\"ahler-Einstein. Moreover, since the automorphism group of $P_1$ is generated by a reflection, by Theorem \ref{tri-auto}, $P_1$ is not symmetric. It has $3$ singular points   of type $A_4+\frac{1}{5}(1,2)+\frac{1}{5}(1,2)$. 
        \item Take $P_2 = conv\{ (11,-3), (0,1), (-11,2) \}.$ Then $P_2$ is K\"ahler--Einstein and not symmetric since $Aut(P_2)$ is trivial. It has $3$ singular points of type $\frac{1}{11}(1,3)+\frac{1}{11}(1,5)+\frac{1}{11}(1,7).$  
    \end{enumerate}
\end{example}

It implies the converse of Theorem \ref{tri-auto}.
\begin{proposition}
    Every subgroup of $D_3$ is an automorphism group of a K\"ahler--Einstein Fano triangle.
\end{proposition}

Similarly one can prove the following.

\begin{proposition}\label{D4D6}
    \begin{enumerate}
        \item Every subgroup of $D_4$ is an automorphism group of a K\"ahler--Einstein Fano polygon with $4$ vertices.
        \item Every subgroup of $D_6$ is an automorphism group of a K\"ahler--Einstein Fano polygon with $6$ vertices.
    \end{enumerate}
\end{proposition}

The proof is given by Examples \ref{KEnotsymm} and \ref{Ex4}.

\begin{example}\label{Ex4}
    \begin{enumerate}
        \item $Aut(P_{0,0}) \cong D_4.$ 
        \item Take a smooth Fano polygon $P=conv \{ (1,0), (1,1),(0,1),$ $(-1,0),$ $(-1,-1),$ $(0,-1) \}.$ Then $Aut(P) \cong D_6$.
        \item Take a singular Fano polygon $P=conv \{ (2,-1),$ $(1,2),$ $(-2,1),$ $(-1,-2) \}$ with $4$ singular points all of type $\frac{1}{5}(1,2)$. Then $Aut(P) \cong \mathbb{Z}/4\mathbb{Z}.$        
        \item Take a singular Fano polygon $P=conv \{ (2,-3),$ $(3,-1),$ $(1,2),$ $(-2,3),$ $(-3,1),$ $(-1,-2) \}$ with $6$ singular points all of type $\frac{1}{7}(1,2)$. Then $Aut(P) \cong \mathbb{Z}/6\mathbb{Z}.$
    \end{enumerate}
\end{example}
 
Recall that finite subgroups of $GL_2(\mathbb{Z})$ are classified.

\begin{theorem}\cite[Theorem 3.3]{Mackie}
    Let $G$ be a finite subgroup of $GL_2(\mathbb{Z})$. Then $G$ is isomorphic to a subgroup of $D_4$ or $D_6$.
\end{theorem}

Thanks to the above theorem, Proposition \ref{D4D6} implies Theorem \ref{KE_aut}

Finally, we remark that the number of vertices of most symmetric Fano polygons is even. According to \cite{GRDB}, up to index $17$, the number of vertices of all K\"ahler-Einstein Fano polygons is even if it is not a triangle. We still do not know any example of a K\"ahler-Einstein Fano polygon with  $5$ or $7$ vertices. But we found such an example with $9$ vertices, which has index $105$. 

\begin{example}
Let $P$ be a Fano polygon with $9$ vertices 
$$vert(P) = \{ (0,1),(-7,3),(-8,3),(-7,2),(0,-1),(3,-2),(7,-3),(7,-2),(5,-1) \}.$$  
Then $P$ is a  symmetric Fano polygon   since $P$ admits a rotation. Indeed,    $$Aut(P)=\Big\{\left( \begin{array}{ccr} 1 & 0 \\ 0 & 1 \\ \end{array} \right),\left( \begin{array}{ccr} -3 & -7 \\ 1 & 2 \\ \end{array} \right),\left( \begin{array}{ccr} 2 & 7 \\ -1 & -3 \\ \end{array} \right)\Big\}.$$  
The singularity types of $P$ are $[3],[5],[4,2],[3],[5],[4,2],[3],[5],[4,2]$ and the index of $P$ is $105$.
\end{example}

\subsection{Demazure root systems}
 
\begin{definition}
Let $P$ be a Fano polytope and $\Delta$ be the corresponding fan in $ N_\mathbb{R}$. Then the set 
    $$\mathcal{R}(P) := \{ r \in M  \; | \; \exists \tau \in \Delta(1)  : \langle v_\tau, r \rangle = -1, \langle v_{\tau'}, r \rangle \geq 0 \; \forall \tau' \in \Delta(1) \setminus \{ \tau \} \}$$
    is called the \emph{Demazure root system} of the Fano polytope $P$. We say that the fan $\Delta$ is said to be {\emph semisimple} if for every $r \in \mathcal{R}$ we have $-r \in \mathcal{R}$.
\end{definition}

\begin{example}\label{Pmn-semisimple}
    The Demazure root system of $P_{m,n}$ is $\{ (1,0), (0,1), (-1,0),(0,-1) \}$ if $m=n=0$, and empty if otherwise. Hence the corresponding fan $\Delta$ is semisimple.
\end{example} 

The Demazure root system is useful in the study of the  automorphism group of the complete toric variety.

\begin{proposition}\cite[Proposition 3.2]{Nill2006}\label{reductive<=>semisimple}
The automorphism group of a complete toric variety is reductive if and only if the corresponding fan is semisimple. 
\end{proposition}

\begin{corollary}
    Let $P$ be a lattice symmetric Fano polygon. Then the automorphism group of the corresponding toric del Pezzo surface is reductive.
\end{corollary}

\begin{proof}
    If $P \ncong P_{m.n}$, then it is symmetric by Theorem \ref{Aut(L-symm)}, so it is K\"ahler--Einstein by Theorem \ref{symm=>KE}, hence it has reductive automorphism group by \cite[Theorem 1.3]{ABHLX}. The fan of $P_{m,n}$ is semisimple by Example \ref{Pmn-semisimple}, hence  the result follows from  Proposition \ref{reductive<=>semisimple}.
\end{proof}

\bigskip {\bf Acknowledgements.}   The authors would like to express their gratitude to the anonymous referees for pointing out errors in previous versions and providing invaluable suggestions that have significantly enhanced the readability of the manuscript, including those leading to the formulation of Proposition \ref{symmetric-dual} and Theorem \ref{symm=>KE2}. D. Hwang would like to thank Grzegorz Kapustka and Micha\l\ Kapustka for informing him the paper \cite{BB} and brought this subject to his attention, Thomas Hall for his interest in this work, and Kyeong-Dong Park for some useful discussions.   
D. Hwang was supported by the Samsung Science and Technology Foundation under Project SSTF-BA1602-03, the National Research Foundation of Korea(NRF) grant funded by the Korea government(MSIT) (2021R1A2C1093787) and the Institute for Basic Science (IBS-R032-D1).

\bigskip

\end{document}